\documentclass[10pt]{amsart}

\usepackage{fullpage}
\usepackage{amsmath, amssymb, amsthm, amsfonts, mathrsfs}
\usepackage{cite}
\usepackage{amscd}
\usepackage{url}
\usepackage{graphicx}
\usepackage{caption}

\usepackage{lscape, pdflscape}

\title{Detecting Concepts Crucial for Success in Mathematics Courses from Knowledge State-based Placement Data}
\author{Marc Harper and Alison Ahlgren Reddy}

\makeatletter
\g@addto@macro{\endabstract}{\@setabstract}
\newcommand{\authorfootnotes}{\renewcommand\thefootnote{\@fnsymbol\c@footnote}}%
\makeatother
\begin{document}
\date{}

\begin{center}
  \LARGE 
  Detecting Concepts Crucial for Success in Mathematics Courses from Knowledge State-based Placement Data \par \bigskip

  \normalsize
  \authorfootnotes
  Marc Harper\footnote{corresponding author, email: marcharper@ucla.edu}\textsuperscript{1} Alison Ahlgren Reddy\textsuperscript{2}\par \bigskip

  \textsuperscript{1}Department of Genomics and Proteomics, UCLA, \par
  \textsuperscript{2}Department of Mathematics, University of Illinois\par \bigskip
  
\end{center}


\begin{abstract}
We show that individual topics and skills can have a dramatic effect on the outcomes of students in various mathematics courses at the University of Illinois. Data from the placement program at Illinois associates a knowledge state, a subset of 182 items and skills that a student is able to complete successfully and repeatedly, with their final grades in a variety of courses from college algebra through multivariate calculus. Using various conditional probabilities and odds ratios, we classify items based on their association with successful and unsuccessful course outcomes, showing that some skills that are advanced for some courses are fundamental or basic to more advanced courses. We examine the impact of specific items across the courses in the traditional college algebra, precalculus, and calculus sequence, as well as courses not typically covered by placement programs, such as higher calculus courses. Visualizations of the knowledge states associated to each student are given for some specific courses and for the entire collection of placement courses, allowing the impact of specific topics to be examined across the undergraduate curriculum.
\end{abstract}


\section{Introduction}

Traditional course placement programs at large universities typically utilize a paper or online exam consisting of approximately 20 to 30 multiple-choice questions. Recently many institutions have begun using assessments based on the theory of knowledge spaces for course placement; at the time of writing, this includes the University of Illinois, the University of Colorado, the University of Texas, Purdue University, the University of Florida, the University of Missouri, and many other higher-education institutions in the United States. These placement programs have produced very interesting datasets. We explore a large such dataset from the University of Illinois.

The theory of knowledge spaces has been developed over the previous two decades and has a major commercial implementation by ALEKS (\textbf{A}ssessment and \textbf{Le}arning in \textbf{K}nowledge \textbf{S}paces). Loosely, a knowledge space is a set of items (essentially problem types) that are related in a hierarchical structure that captures the dependence between items. Students can be assessed in an adaptive manner to determine the precise subset of the total domain of items that the student is able to successfully complete. This subset is called a \emph{knowledge state}. Not all subsets are feasible knowledge states, and the ability to complete any particular item has a probabilistic relationship with the ability to complete any other item, which can be determined empirically. The commercial implementation of knowledge spaces by ALEKS can assess students on domains of 200 to 300 items with just 30 questions. For information on the theory of knowledge spaces, see \cite{falmagne2011learning} \cite{doignon1994knowledge}; for the practical usage of the theory see \cite{falmagne2007assessing} \cite{falmagne2006assessment} \cite{falmagne1990introduction}.

The typical implementation of a knowledge space-based placement program is for each student to take an assessment before the beginning of his or her next semester or quarter, place into and enroll in a course, and ultimately pass or fail the course. Essentially every student has a unique knowledge state since it is highly unlikely that any two students in a particular course have the exact same knowledge state (unless it is empty or the complete domain). At the time of writing, for practical reasons all institutions using knowledge states for course placement reduce the state to a score by taking the cardinality of the state and dividing by the total number of items in the domain, assigning to each the percentage of the domain that he or she demonstrated mastery of in the assessment. Having a knowledge state for each student rather than simply a score for a 20-30 question paper-and-pencil test provides a wealth of additional information. In this manuscript we investigate what can be learned from the content of students knowledge states rather than just their size.  Some analyses on subsets and summaries of this data (based on the scores of the knowledge states) have been previously published \cite{ahlgren2013mathematics} \cite{reddy2013aleks}, and this is the first to look at the statistics of individual topics.

In particular, we describe statistics that are useful for detecting items that have high impact on student outcomes. More precisely, we identify specific items that are highly correlated with success (passing with a grade of C- or better) or failure in various courses, which is of interest to instructors, course planners, book writers, standards writers, and mathematics education researchers, and educational policy makers. Since all incoming students take the same assessment, we can compare the knowledge states of students in different courses, which is often not possible with other placement mechanisms. The authors hope that this work leads to a closer investigation of data at other institutions and ultimately to large cross-institutional studies that may determine which items are universally-linked to course outcomes, and caution that the results presented in this manuscript are specific to the student population and courses at the University of Illinois.

\section{Description and Visualization of Data}

The placement mechanism at the University of Illinois is an ALEKS assessment on a customized version of the ALEKS Preparation for Calculus domain, which is an adaptive assessment on a collection of 182 items (or problem types).\footnote{Please note that there are several ALEKS domains in use as placement assessments, including the full Preparation for Calculus domain with over 200 items, and a recently released placement-targeted domain with over 300 items, intended for placement over a larger range of courses. The complete list of items in the domain Preparation for Calculus domain, including explicit problem instances, is available at \url{http://www.aleks.com/}.} This particular domain has items ranging from basic algebra through precalculus, but no actual calculus content (no derivatives, for instance). The end result of an adaptive assessment is a knowledge state, which is simply the subset of the 182 items for which the student demonstrated mastery, as determined by the assessment inference mechanism. Each time a student responds to a question, the probability associated to every other item in the domain is adjusted, and the assessment ends after approximately 30 questions. Our data set consists of an initial knowledge state within four months of course start for each student and each class in the placement program (and a few others) since 2008, as well as the final grade achieved by each student at the end of the course. We will only discuss fall semesters because they typically have larger course sizes. For more information on the placement program at the University of Illinois see \cite{ahlgren2011assessment}.

Students are placed into a course based on the total percentage of items in their knowledge state, which we will refer to as the score. Although only five courses (College Algebra, PreCalculus, Business Calculus, Calculus I, and Calculus I for experienced students) are governed by the placement program, all incoming students are instructed to complete at least one ALEKS assessment after acceptance to the university, and most students in courses not governed by the placement program complete an assessment. This means that for some courses outside of the placement program, including Calculus II, Calculus III, and Elementary Linear Algebra (for non-technical students), most students have completed initial assessments that can also be used for analysis.\footnote{A description of these courses is available online at \url{http://www.math.illinois.edu}.}  In particular, we note that for the placement courses, students have been filtered by a minimum score on the placement assessment, but for courses outside the placement program, there are no such restrictions on student enrollment.

\begin{figure}
    \centering
    \includegraphics[width=0.4\textwidth]{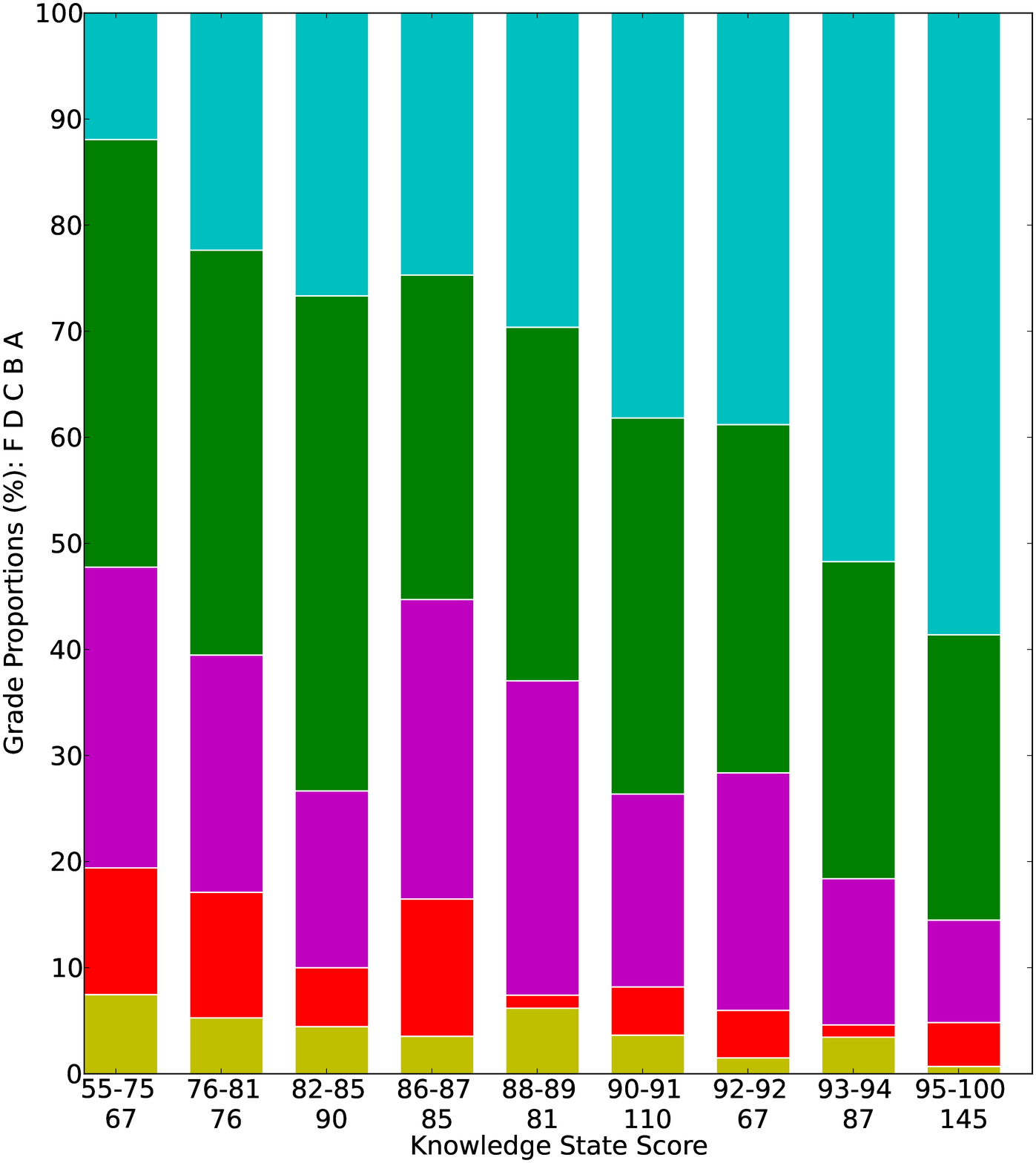}
    \includegraphics[width=0.4\textwidth]{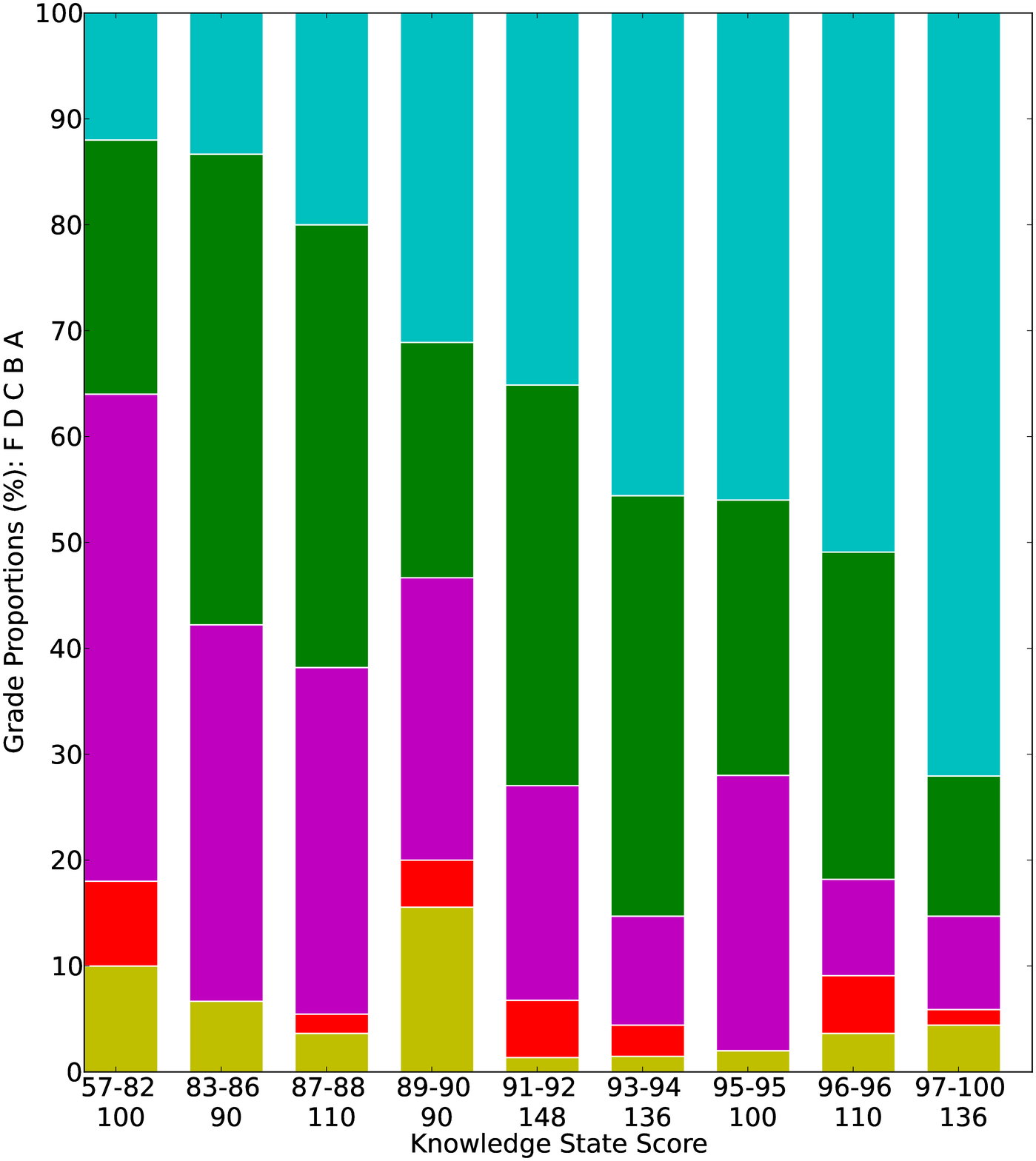}
    \caption{ALEKS Score (knowledge state percentage of total) versus grades. Left: Calculus II (231). Right: Calculus III (241, multivariate). From top to bottom, The colored bars indicate the proportion of A/A-, B+/B/B-, C+/C/C-, D+/D/D-, and F/W grades for each grouping. The numbers under each group indicate the range of scores and number of students in the group. Both courses are from the Fall semester of 2010, and in neither have any placement requirements. Notice that the proportion of A and B grades increase dramatically as the placement scores increase even though the assessment covers no content specific to calculus (such as derivatives and integrals). See \cite{ahlgren2013mathematics} for similar figures for placement courses (college algebra through calculus I) and linear algebra.}
    \label{calc_2_3}
\end{figure}

In Figure \ref{calc_2_3}, we show the aggregate distribution of ALEKS scores in two courses not governed by the placement program: Calculus II and Calculus III (multivariate) for the year 2010. Perhaps surprisingly, despite the fact that there is no explicit calculus content in the assessment, there is a strong and obvious relationship between the initial ALEKS score and the distribution of grades of students in the course. This suggests that success in subsequent calculus courses is partially based on what are generally considered non-calculus or arithmetic and algebra skills rather than on basic calculus knowledge. This phenomenon was previously observed in a small and informal study by Stephen Wilson \cite{wilson}. Similar plots are available for several other courses and years in \cite{ahlgren2013mathematics}, and the pattern shown in Figure \ref{calc_2_3} is essentially always present, and often more dramatically so. 

\begin{landscape}

\begin{figure}[h]
\centering 
\includegraphics[width=1.4\textwidth]{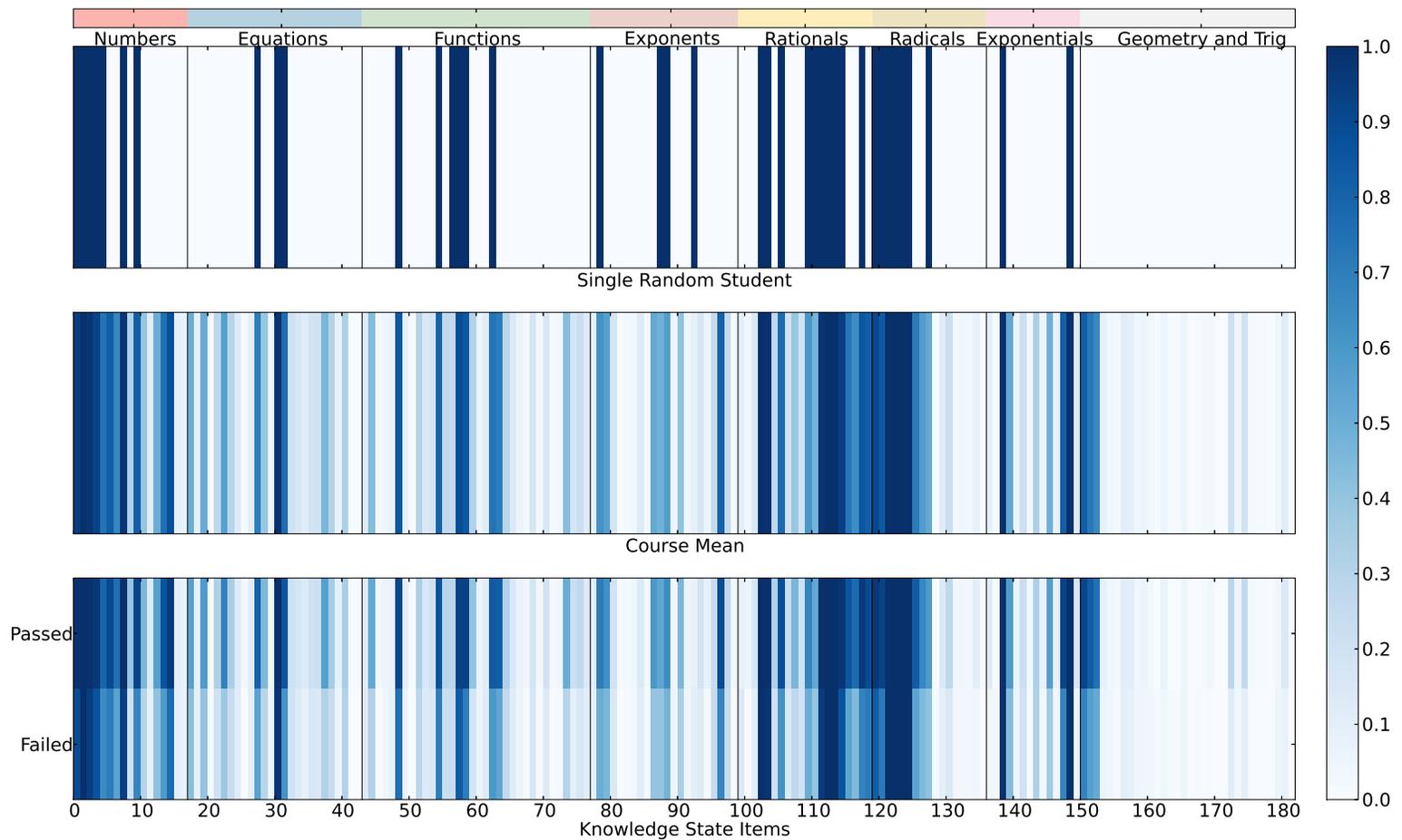}
\caption{Visualization of knowledge states in College Algebra (227 students). Top: One randomly chosen student, having approximately 40\% of the items in the his or her knowledge state. Middle: Average knowledge state, color indicates proportion of students having the item in their knowledge state. There are many items (white) that no student has, since the set of items goes significantly beyond the course material and expected initial knowledge for College Algebra. Similarly, some items are significantly below the course level, and all or nearly all students have the item (dark blue). Bottom: Average knowledge states for students split into passing (C- or better) and failing. Some but not all items skew in proportion among the passing and failing students.}
\label{barcode_pass_fail}
\end{figure}

\end{landscape}

\begin{landscape}

\begin{figure}[h!]
\includegraphics[width=1.4\textwidth]{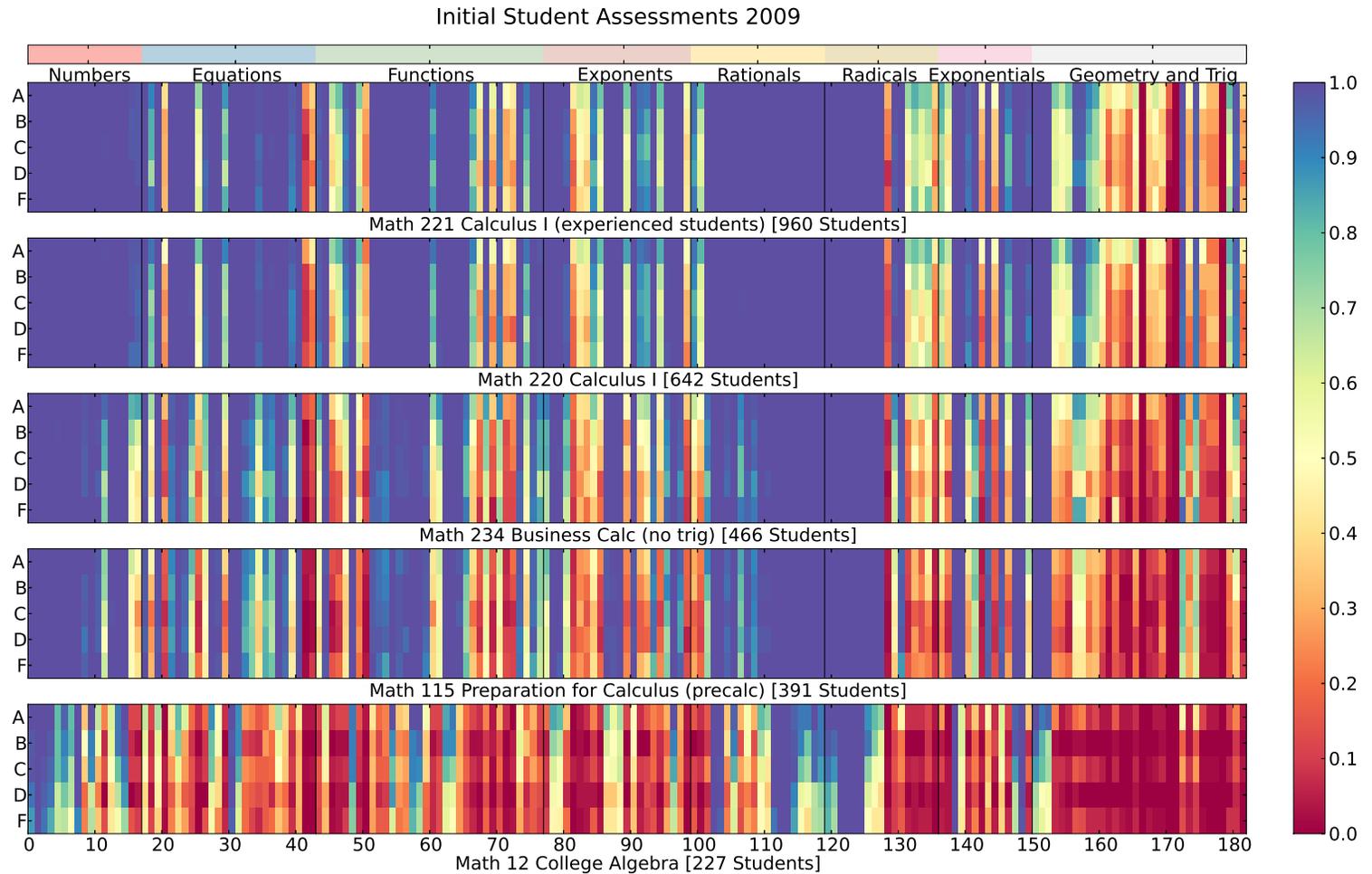}
\caption{Visualization of Knowledge States from College Algebra through Calculus I for experienced students in the Fall semester of 2009. Color indicates proportion of students with the item in their knowledge state for the particular grade. Many items vary over one or more classes, and many items that vary in College Algebra and PreCalculus have little impact in Calculus. Some items slowly increase in proportion over all classes from Bottom to Top. Courses, from Bottom to top: College Algebra (012), Preparation for Calculus / PreCalculus (115), Business Calculus (234), Calculus I (220), and Calculus I for experienced students (221).}
\label{barcode_grade_2009}
\end{figure}

\end{landscape}

For each student in each course in the fall semesters of the years 2008-2011, we associate the knowledge state of the best ALEKS assessment before the course and the final grade in the course, given as a letter $A$, $A-$, $B+$, etc. Passing is defined as a $C-$ or better in all courses. Each knowledge space is encoded as a binary vector of ones and zeros corresponding to whether or not the knowledge state includes each item. Enumerating the items in the knowledge space from 0 to 181, we can visualize a single knowledge space as a \emph{barcode}, with the bar filled in if the student has the item and empty otherwise. From these vectors we can compute the proportion of students in each course having a particular item, and split the students into those that passed and those that failed, and compare the mean barcodes. See Figure \ref{barcode_pass_fail}. Notice that some items vary significantly between the two groups whereas others do not. One reason that this can occur is that for a given course, some of the items are essentially not in play since either all students in the course have the item in their knowledge states (the course is beyond the particular item) or no students have the item (the item is an advanced topic that students in the course have not yet encountered). Note that while the visual representation suggests that items are independent, they are generally not. Figure \ref{barcode_grade_2009} shows students in the sequence of placement courses, in order of perceived difficulty from College Algebra through Calc I, split by grade instead of passing and failing. Since all students take an assessments on the same set of topics, we can clearly see how particular items vary over the initial knowledge of the various student populations from course to course.

\begin{landscape}

\begin{figure}[h!]
\includegraphics[width=1.4\textwidth]{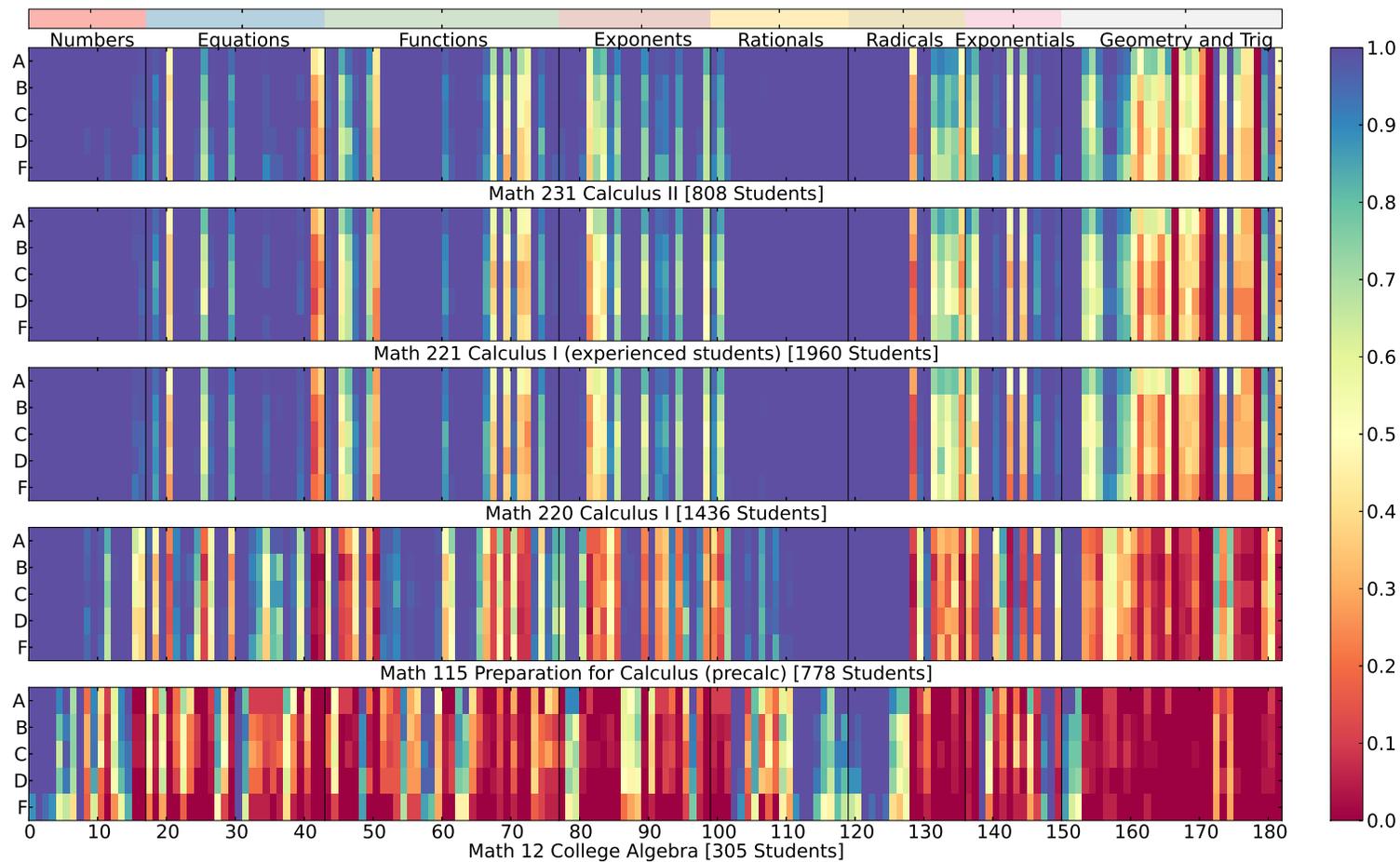}
\caption{Visualization of Knowledge States from College Algebra through Calculus II for the Fall semester of 2010. Compare to Figure \ref{barcode_grade_2009} (note that Business Calc was not offered in Fall 2010). As before, color indicates proportion of students with the item in their knowledge state for the particular grade. Courses, from Bottom to top: College Algebra (012), Preparation for Calculus / PreCalculus (115), Calculus I (220), and Calculus I for experienced students (221), Calculus II (231).}
\label{barcode_grade_2010}
\end{figure}

\end{landscape}

Consider Figures \ref{barcode_grade_2009} and \ref{barcode_grade_2010} in detail. It is easy to see that some items are known by every student in some courses, denoted by a vertical strips of all purple, while other items are known by no students are denoted a vertical strips of red. Many items vary with grade not just over a single course, such as items 114 through 121 in College Algebra, but over multiple courses and even the entire curricula, such as several items numbered between 40 and 50 in the functions category which include problem types on topics such as domain and range. Students in College Algebra are unlikely to have these items in their knowledge state at the beginning of the course, but for calculus students, those earning higher grades are more likely to have these items in their knowledge states. Item 18, solving a one-variable linear equation, has a similar behavior.


\section{Methods and Results}

Results are qualitatively similar for each year examined, and we will mention data from 2008 and 2011 when it is particularly interesting, but otherwise omit these years to avoid redundant diagrams and statistics. We note that absolutely no attempt was made to normalize grading rubrics or scales or to have consistent instructors, exams, or course materials across the courses and years in the study, so the data comes from the placement program \emph{in natura}. Finally, note that in our data class sizes increase as the course number / difficulty increases: College Algebra has just under 200 students in contrast to over 1000 students in some of the Calculus courses, so the data from the lower level courses tends to be noisier. Because there are many items and many courses, we develop methods to systematically detect interesting items without resorting to manual examination, and to classify items.

\subsection{Odds Ratios}

\begin{figure}[h!]
\centering 
\includegraphics[width=0.6\textwidth]{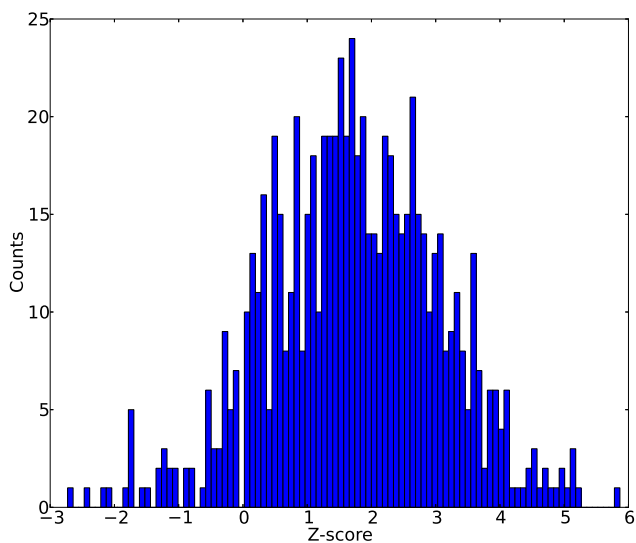}
\caption{Z-scores for log odds ratios of every item in every class (each item potentially appears multiple times) in the Fall of 2010 (all classes mentioned above except Business Calc, which was moved to the spring semester in 2010). Items that either all students in a course have or all students in a course do not have are excluded, as well as any items with $n_{00}=0$ or $n_{11}=0$. 91\% of the items have log odds ratio greater than zero, with mean 1.74 and standard deviation 1.3 over all classes. Other years are qualitatively similar. The standard deviation drops to $~\approx 0.9$ if small counts items are excluded (less than 30 students having or not having the item).}
\label{zscores_all}
\end{figure}

A basic statistical tool for categorical data is the odds ratio \cite{agresti2007introduction} \cite{rao2006handbook}. We calculate the odds ratio from contingency tables (Figure \ref{contingency_def}) for each particular item. The odds ratio is defined as
\[ OR = \frac{n_{00} n_{11}}{n_{01} n_{10}}\]
A ratio greater than one indicates a positive correlation between the two binary variables of having an item (or not) in the student's knowledge state and course performance (passing the course or not). An odds ratio of precisely one would indicate no significant relationship between the variables. For small counts, it is possible for $n_{01}$ or $n_{10}$ to be zero. For this reason we add a constant $C=0.5$ to all four entries to create pseudocounts. This does not alter the asymptotics of the estimate and prevents undesirable values of the odds ratio while still giving reasonable estimates \cite{agresti2007introduction}. The log odds ratio is distributed normally with mean $\log OR$ and variance given by
\[ \sigma^2 = \frac{1}{n_{00}} + \frac{1}{n_{01}} + \frac{1}{n_{10}} + \frac{1}{n_{11}}. \]
We can empirically estimate the odds ratio and use it to detect interesting items in a given course. For any given item and course, the log odds ratio would be zero if possession of the item and passing the course were independent. Figure \ref{zscores_all} shows that the log odds ratios are not distributed about zero as would be expected in this neutral case, and that many items have large associations with course performance. Indeed, aggregated over all classes in 2010, approximately 91\% of the odds ratios of items with $n_{00}, n_{11} \neq 0$ are greater than 1. In light of Figure \ref{calc_2_3} this is perhaps not surprising, but it is not a foregone conclusion that most items correlate positively with course performance just because the aggregate total over all items correlates, since a small subset of items could be vastly more informative than others.

\begin{figure}[h!]
\centering 
\begin{tabular}{r|c|c}
    & Passed & Failed\\
    \hline
    has item & $n_{00}$ & $n_{01}$ \\
    \hline
    not has item   & $n_{10}$ & $n_{11}$ \\
    \hline  
\end{tabular}
\caption{Contingency table of counts of coincidences of students passing or failing and having or not having a particular item in the knowledge space. Such tables were computed for every item in every course and used to calculated the odds ratio and various other quantities.}
\label{contingency_def}
\end{figure}

\subsection{Two Classifying Properties for Evaluating Items}
In addition to the odds ratio, which we can use to systematically identify strongly correlated items, we can also calculate various useful probabilities. From the contingency tables, we can compute the following two ratios of conditional probabilities:
\begin{align*}
&\text{Passing Ratio:}&  &\frac{\text{Pr}(\text{Pass} | \text{has item})}{\text{Pr}(\text{Pass} | \text{not has item})} \\
&\text{Failing Ratio:}&  &\frac{\text{Pr}(\text{Fail} | \text{not has item})}{\text{Pr}(\text{Fail} | \text{has item})}
\end{align*}

If the passing ratio is greater than one, then having the item is associated with the student passing the course. We will refer to an item as having the \emph{passing property} if the ratio is greater than one. Similarly, if the failing ratio is greater than one, we refer to the item as having the \emph{failing property}. Not having such items are associated with a greater chance of failure. The product of these ratios is the odds ratio:
\[ OR = \frac{\text{Pr}(\text{Pass} | \text{has item})}{\text{Pr}(\text{Pass} | \text{not has item})} \frac{\text{Pr}(\text{Fail} | \text{not has item})}{\text{Pr}(\text{Fail} | \text{has item})}\]
The size of the ratio indicates how strongly the item skews with passing and/or failing. So a large odds ratio could indicate either a strong relationship with having an item and passing, with lacking an item and failing, or both. We can further classify items by whether or not they have one or both of the two properties, and how strongly so.
\begin{itemize}
    \item \emph{Basic:} Basic items are those that have large failing ratios. Students lacking these items are significantly more likely to fail yet students having the item are not particularly more likely to pass. This indicates that the item is basic or foundational to the course material.
    \item \emph{Advanced:} Advanced items are those that have large passing ratios. Students with these items are significantly more likely to pass but are not significantly more likely to fail if they lack the item. This indicates that the item corresponds to a more advanced topic.
    \item \emph{Core:} Items with both properties are called core items. These items detect both well-prepared students and poorly-prepared students, but typically to a lesser extent in either case than a basic or advanced item.
    \item \emph{Uninformative:} Items with neither property are called uninformative; this includes items which every student in the course has or does not have, which give no information regarding course performance. In this degenerate case the ratios are ill-defined, either zero or infinity, or there is not enough data to give a good estimate.
\end{itemize}

An item that is uninformative for one course need not be so for another course. Such an item might be very basic or very advanced in comparison to all items in the knowledge space, but relatively inconsequential for the course under consideration, since there is little if any variation in possession of the item among students in the course. In our data, most items that have one of the two properties have the other, but not necessarily to the same extent. It is useful to compare the conditional probabilities to the ambient rates, e.g. $\text{Pr}(\text{Pass} | \text{has item})$ to the proportion of students that passed the course unconditionally.

Figure \ref{item_classifications} shows the classification of items by course. Interesting, many items with the passing property also have the failing property, and it appears that the failing property is the greater contributor to the odds ratio. This is simply because having a particular item does not generally confer a particularly strong benefit to course performance, but lacking an item fundamental to the course can have a substantially negative effect on course outcomes. On the other hand, as we remarked above in Figure \ref{calc_2_3}, the aggregate total of passing items in a knowledge state is strongly associated with course outcomes. More plainly, a well-prepared student is indicated by a breadth of knowledge covering the essential topics and subcategories (such as trigonometry and logarithms), whereas a poorly-prepared student is indicated by the lack of just one or a small number of basic items (even if their total number of items is relatively large), and this is reflected in the magnitude of the failing and passing ratios for each item.

\begin{figure}
\centering
\begin{tabular}{|c|c|c|c|c|c|c|}
\hline
Course & Year & Passing & Failing & Both & Neither & Uninformative\\
\hline \hline
College Algebra & 2009 & 146 & 139 & 136 & 32 & 13\\ \hline
 & 2010 & 134 & 110 & 104 & 41 & 33\\ \hline
PreCalculus & 2009 & 81 & 82 & 80 & 98 & 54\\ \hline
 & 2010 & 89 & 89 & 88 & 91 & 52\\ \hline
Business Calculus & 2009 & 113 & 114 & 113 & 67 & 54\\ \hline
Elementary Linear Algebra & 2010 & 105 & 108 & 105 & 73 & 58\\ \hline
Calculus I & 2009 & 76 & 77 & 76 & 104 & 92\\ \hline
 & 2010 & 80 & 82 & 80 & 99 & 96\\ \hline
Calculus I Exp* & 2009 & 77 & 74 & 74 & 104 & 91\\ \hline
 & 2010 & 69 & 67 & 66 & 111 & 91\\ \hline
Calculus II & 2010 & 87 & 88 & 84 & 90 & 81\\ \hline
Calculus III & 2010 & 70 & 69 & 67 & 109 & 94\\ \hline
\end{tabular}
\caption{Item classifications for various courses. Of the 182 items in the assessment, the number of items with the Passing, Failing, or both properties are listed. The remaining items are separated into those that have neither property, and those that are uninformative. *Students with prior calculus experience can self-enroll in Calculus I Exp rather than Calculus I.}
\label{item_classifications}
\end{figure}

\subsection{Examples}

Consider the following examples of items fitting the above classification for multiple years for the same course in our dataset.

\begin{enumerate}
\item[(alge026)] \emph{College Algebra:} Students are asked to simplify a quotient of products of powers, such as 
\[ \frac{wu^6}{w^4u^2} \]
In the Fall of 2010, 88\% of the incoming students were able to complete this problem and passed 75\% of the time (the ambient pass rate was 71\%). Of the 12\% of students unable to complete this item, 57\% failed the course. The odds-ratio z-score is 3.87. This is a mostly basic item, detecting very unprepared students, and had similar statistics in 2009.
\item[(alge032)] \emph{College Algebra:} Students are asked to expand a binomial, such as $(5t+6)^2$.\\ This item has similar statistics to the previous item with an odds-ratio z-score of 4.6 in 2009 and 3.7 in 2010. Students lacking this item were at least twice as likely to fail.
\item[(geom044)] \emph{All courses} In 2010, an item asking for students to compute the hypotenuse of a right triangle given the other sides had an odds-ratio z-score of 3.1 in Precalculus, i.e. to apply the Pythogorean theorem. This is item \#165 in Figures \ref{barcode_grade_2009} and \ref{barcode_grade_2010}, and it behaves similarly for Business Calculus and the two flavors of Calculus I. It is a basic or core item for these courses. For College Algebra in 2010, this item had z-score of 4.0 and is an advanced item as relatively few students have it (27\%) in their initial knowledge state but 90\% of them pass!
\item[(fun030)] \emph{Business Calculus / Calculus:} Students are asked to compute the outputs of a piecewise-defined function at specific points (e.g. $g(-1)$ and $g(2)$) for a piecewise polynomial function. In Business Calculus in 2008, 24\% of the students had this item in their knowledge state and passed the course 94\% of the time versus 80\% overall. Of the students lacking the item 24\% failed the course (versus 20\% overall). This is an advanced item for this course indicating well-prepared students and the item has similar statistics in 2009. This item has large z-score in several Calculus courses and years (even Calculus II and III).
\item[(alge712)] \emph{Calculus I, experienced:} Students are asked to graph an exponential function, e.g. 
\[ g(x) = -\left(\frac{3}{5}\right)^x\]
In the Fall of 2008, about half of the students have this item in their knowledge states; of those, 91\% passed the course (versus 85\% overall) and of the students lacking the item, 21\% fail (versus 15\%) overall. This is a core item that has both the passing and failing properties, and has replicate statistics in 2009 and 2010 (z-score 4.8 in 2008).
\item[(alge682)] \emph{Calculus II:} Students are asked to simplify a ratio of polynomials, e.g.
\[ \frac{x^2 -9x +20}{48-3x^2}\]
In both 2010 and 2011, approximately 90\% of the students passed this course, and 95\% of the students had this item. Of the 5\% lacking the item, 20\% failed in 2010 and 35\% failed in 2011!
\end{enumerate}

In the last example for Calculus II, it is not surprising that this item is a basic item and that lack of the item is strongly correlated with failing the course. This course typically covers several topics in which the ability to manipulate rational functions is critical, such as integration of rational functions, particular the method of partial fractions. This item also has large z-scores for Business Calculus, which may be due to the reliance on rational functions for problems and examples as a substitute for trigonometry-based examples. On the other hand, that the understanding of piecewise functions is strongly correlated with passing Business Calculus is less intuitively clear since it is unlikely to consume a large amount of class time. A familiar with piecewise functions may indicate adeptness with functions in general, certainly important for Calculus of any variety. Similarly, the ability to graph an exponential function may not be a crucial for success skill per se, yet may indicate an understanding of functions and their behaviors that is reflected in the ability to compute limits and other critical skills in calculus that are based on the shape and geometry of functions.

There are many other items with interesting statistics and it is not our intention to speculate on the reasons for the properties of every such item, so we focus on a few common outcomes. One particularly common scenario seen above is an item that most students in a course have in their knowledge state but a small percentage (less than 10\%) do not have. For instance, in one calculus course (Fall 2010), 3\% of the students lacked an item that simply asked them to compute $\log_3{27}$ (or numerically similar variant), and these students were more than twice as likely to fail the course. Commonly there are several such items in any particular course, and they seem to be at the base of a collection of items forming a significant ``branch'' of the domain (in this case, the item is one of the more simple items in the exponentials and logarithms subcategory). 

This highlights an omission in the basic preparedness in an important topic area, and suggests that students may fail more advanced courses like calculus not because of the content of the course itself, but because of a lack of ability to manipulate the underlying objects to which calculus is applied. This suggests that a closer examination of a student's knowledge state beyond the score may improve placement effectiveness, and could be utilized by academic advisors and instructors, for instance, by enabling detailed conversations with students regarding the likelihood of their success in various courses.

\section{Discussion}

In this manuscript we have shown that it is possible to connect student outcomes in standard mathematics courses with individual items from the initial knowledge of students. By identifying items that are highly correlated with student outcomes, curricula adjustments and other policy decisions can be data driven. For instance if a particular item is basic for Calculus I but students are not entering from precalculus with this item mastered, the precalculus curriculum can be adjusted to place additional emphasis on the item.

Broadly speaking, our analysis strongly suggests that knowledge state-based placements provide a substantially greater amount of information about student preparedness and course performance than traditional placement exams, such as fixed question paper-and-pencil tests (and their online variants). Such exams are typically no more than 30-40 questions (and often a few as 20), and so cannot possibly provide the granularity of an adaptive assessment over a range of courses. Our results confirm anecdotal observations that students in courses as high as multivariate calculus can lack very basic skills and are at substantially higher risk of failing because of their underpreparedness.

The authors strongly caution against generalizing from the results presented here until corroborating results are found from other institutions. Until such studies are performed, it can only be said that the items explicitly presented here are significant to the student body and courses at the University of Illinois. Future work may elucidate items that are universally critical for success in standard courses and the depth of local variation of student preparation and instruction at diverse institutions.

\section{Acknowledgements}

\subsection{Thanks} The authors are indebted to Alison Champion of the University of Illinois for assistance in data collection.

\subsection{Figures} All plots created with \emph{matplotlib} \cite{Hunter:2007}.

\subsection{IRB} Student data was collected under the guidelines of Institution Review Board \#07752.

\subsection{Conflicts of Interest} The authors have been compensated consultants for the implementation of placement programs at various institutions. Marc Harper consults at ALEKS on the development of new mathematics placement software.

\bibliography{ref}
\bibliographystyle{plain}

\end{document}